\documentclass[12pt]{article}
\usepackage{amsmath}
	\allowdisplaybreaks
\usepackage{amsthm,amssymb,mathrsfs}
\usepackage{stmaryrd}
\usepackage{dsfont}
\usepackage{mathtools}
\usepackage[shortlabels]{enumitem}
	\setlist[enumerate]{resume}
	\setlist[enumerate,itemize]{topsep=0em,itemsep=.5em,partopsep=0em,parsep=0em}	
\usepackage{graphics}
\usepackage[margin=3cm]{geometry}
\usepackage{changepage}
\usepackage{multicol}

\usepackage{color,soul}

\usepackage{manfnt}
\reversemarginpar

\usepackage[colorlinks=true,
    bookmarksnumbered=true, linktocpage=true,
	 linkcolor = blue, citecolor = magenta, urlcolor = black]{hyperref}
\usepackage{bookmark}

\usepackage{tikz}

\usepackage{tikz-cd}
	\tikzcdset{every label/.append style = {font = \small}}

\usepackage{titlesec}
	\titleformat{\section}{\bf\filcenter}{\thesection.}{.5em}{}
	\titleformat{\subsection}{\bf}{\thesubsection}{0em}{}
	\titlespacing*{\section}{0pt}{1em}{.5em}
	\titlespacing*{\subsection}{0pt}{.5em}{.25em}

\theoremstyle{definition}
\newtheorem{definition}{Definition}[]  
\newtheorem{context}[definition]{\indent}

\numberwithin{equation}{definition}

\title{\normalsize\bf ON ENDOMORPHISM ALGEBRAS OF GELFAND--GRAEV REPRESENTATIONS II\footnote{This article has been published on the Bull.\;London Math.\;Soc.\;(\url{https://doi.org/10.1112/blms.12899}) under the licence \href{https://creativecommons.org/licenses/by/4.0/legalcode}{CC BY 4.0}. The present version is the accepted version for the BLMS, with the section numbering changed in accordance with the published version on the BLMS.}}
\author{\normalsize Tzu-Jan Li and Jack Shotton}
\date{}
\begin{document}

\maketitle

\vspace{-15mm}
\begin{abstract}
  \noindent{\bf Abstract.} Let $G$ be a connected reductive group defined over a finite field $\mathbb{F}_q$ of
  characteristic $p$, with Deligne--Lusztig dual $G^\ast$. We show that, over $\overline{\mathbb{Z}}[1/pM]$ where $M$ is the
  product of all bad primes for $G$, the endomorphism ring of a Gelfand--Graev representation of $G(\mathbb{F}_q)$ is
  isomorphic to the Grothendieck ring of the category of finite-dimensional $\overline{\mathbb{F}}_q$-representations
  of $G^\ast(\mathbb{F}_q)$.
\end{abstract}

\begin{context}\label{section-intro}

  {\bf Introduction.}\;--- Let $G$ be a connected reductive group defined over a finite field $\mathbb{F}_q$ of
  characteristic $p$, let $F$ be the associated Frobenius endomorphism of $G$, and let $\Lambda$ be a subring of
  $\overline{\mathbb{Q}}$ containing $\overline{\mathbb{Z}}[\frac{1}{p}]$. Let $B_0$ be an $F$-stable Borel subgroup of
  $G$ with (necessarily $F$-stable) unipotent radical $U_0$, and let $\psi : U_0^F \longrightarrow \Lambda^\times$ be a
  regular (also called nondegenerate) character. The Gelfand--Graev representation
  \[\Gamma_{G,\psi}:=\mathrm{Ind}_{U_0^F}^{G^F} \psi\]
  is an important representation of $G^F$ (already studied in \cite[Sec.\;10]{Deligne--Lusztig} and
  \cite{Digne-Lehrer-Michel}).  Its endomorphism ring
  \[\Lambda\mathsf{E}_G:=\mathrm{End}_{\Lambda G^F}(\Gamma_{G,\psi})\] is commutative, independent of the choice of
  $\psi$ up to isomorphism and, over $\overline{\mathbb{Q}}$, may be identified with the ring of
  $\overline{\mathbb{Q}}$-valued class functions on $G_{\mathrm{ss}}^{\ast F^\ast}$, where $(G^\ast,F^\ast)$ is a chosen
  Deligne--Lusztig dual of $(G,F)$ (see \cite{Curtis}). Such an identification only depends on choices
  of group homomorphisms
  $(\mathbb{Q}/\mathbb{Z})_{p'}\simeq\overline{\mathbb{F}}_q^\times\hookrightarrow\overline{\mathbb{Q}}^\times$, which
  we fix from now on.

  There are then (at least) two natural $\Lambda$-lattices in $\overline{\mathbb{Q}}\mathsf{E}_G$:
  $\Lambda \mathsf{E}_G$ and the lattice $\Lambda \mathsf{K}_{G^\ast}$ spanned by Brauer characters of irreducible
  representations of $G^{\ast F^\ast}$; here, $\mathsf{K}_{G^\ast}$ is the Grothendieck ring of the category of
  finite-dimensional $\overline{\mathbb{F}}_q G^{\ast F^\ast}$-modules. Denoting by $G_{\mathrm{ss}}^{\ast F^\ast}\!/\!\sim$ the set of semisimple conjugacy classes in
  $G^{\ast F^\ast}$, we may then, as in \cite[Sec.\;2.5]{Endomorphism}, identify
  \begin{equation}\label{EK-identify}
    \overline{\mathbb{Q}}\mathsf{E}_G=\overline{\mathbb{Q}}^{G^{\ast F^\ast}_{\mathrm{ss}}\!/\sim}=\overline{\mathbb{Q}}\mathsf{K}_{G^\ast}
  \end{equation}
as $\overline{\mathbb{Q}}$-algebras, where we recall that the second equality follows from the Brauer character
isomorphism $\overline{\mathbb{Q}}\mathsf{K}_{G^\ast}\xrightarrow{\;\sim\;}\overline{\mathbb{Q}}^{G_{p'}^{\ast
    F^\ast}\!/\sim}$ and from the fact that ${G_{p'}^{\ast F^\ast}\!/\!\sim}={G^{\ast F^\ast}_{\mathrm{ss}}\!/\!\sim}$. Here ${G_{p'}^{\ast F^\ast}\!/\sim}$ is the set of $p$-regular conjugacy classes in $G^{\ast F^\ast}$.

  The main result of this paper may now be stated as follows:

  {\bf Main theorem.} {\it If all bad primes for $G$ are invertible in $\Lambda$, then the two $\Lambda$-lattices
    $\Lambda\mathsf{E}_G$ and $\Lambda\mathsf{K}_{G^\ast}$ of $\overline{\mathbb{Q}}\mathsf{E}_G$ are equal.}

  Here, we use the notion of ``bad primes for $G$'' from \cite{Springer66}. Denoting by $R$ the root system of $G$, a
  prime number $\ell$ is called {\it bad} for $G$ if one of the following three conditions holds: (i) $\ell=2$, and $R$
  has an irreducible factor not of type $A$; (ii) $\ell=3$, and $R$ has an irreducible factor of exceptional
  type ($G_2$, $F_4$, $E_6$, $E_7$, or $E_8$); (iii) $\ell=5$, and $R$ has an irreducible factor of type $E_8$.

  In this theorem, the assumption on the bad primes for $G$ is due to the use of {\it almost characters} in Lusztig's
  work on unipotent characters, where bad primes appear in the denominators of the ``Fourier transform matrix.'' We
  expect that the theorem remains true without this assumption, though our present method cannot prove it.

  Our theorem improves the equality
  $\overline{\mathbb{Z}}[\frac{1}{p|W|}]\mathsf{E}_G=\overline{\mathbb{Z}}[\frac{1}{p|W|}]\mathsf{K}_{G^\ast}$ (where
  $W$ is the Weyl group of $G$) in \cite[Thm.\;2.3]{Endomorphism} whenever the adjoint group of $G$ is simple of type
  other than $F_4$ or $G_2$ (in these two excluded types, the bad primes and the primes dividing the order of the Weyl
  group coincide). Moreover, via the $\mathbb{Z}$-model $\mathsf{E}_G$ of $\Lambda\mathsf{E}_G$ from
  \cite[Sec.\;1.5]{Endomorphism}, if we denote by $M$ is the product of all bad primes for $G$, then the above theorem
  implies that $\mathbb{Z}[\frac{1}{pM}]\mathsf{E}_{G}=\mathbb{Z}[\frac{1}{pM}]\mathsf{K}_{G^\ast}$. Indeed, this amounts to showing that the identification $\overline{\mathbb{Z}}[\frac{1}{pM}]\mathsf{E}_{G}=\overline{\mathbb{Z}}[\frac{1}{pM}]\mathsf{K}_{G^\ast}$ in the above theorem is equivariant under the action of the Galois group $\mathrm{Gal}(\overline{\mathbb{Q}}/\mathbb{Q})$ on the coefficients, and the proof of this equivariance is the same as that of \cite[Cor.\;2.4]{Endomorphism}.

  {\bf Relation with invariant theory.} Let $\mathsf{B}_{G^\vee}$ be the ring of functions of the $\mathbb{Z}$-scheme
  $(T^\vee\sslash W)^{F^\vee}$, where $(G^\vee,T^\vee)$ is the split $\mathbb{Z}$-dual of $(G,T)$ with $T$ an $F$-stable
  maximal torus of $G$, $W=N_{G^\vee}(T^\vee)/T^\vee$ is the Weyl group of $(G^\vee,T^\vee)$, and
  $F^\vee:T^\vee\longrightarrow T^\vee$ is induced by the action of $F$ on $Y(T^\vee) = X(T)$.  If $G^\ast$ has
  simply-connected derived subgroup, then $\Lambda\mathsf{B}_{G^\vee}$ is also a $\Lambda$-lattice of
  $\overline{\mathbb{Q}}\mathsf{E}_G$ and appears to be significant for the {\it local Langlands correspondence in
    families}. Indeed, for $\mathrm{GL}_n$, in the course of constructing this correspondence in joint work with Moss
  \cite{Helm-Moss}, Helm proved in \cite[Thm.\;10.1]{Helm} the equality
  $\Lambda\mathsf{E}_{\mathrm{GL}_n} = \Lambda \mathsf{B}_{\mathrm{GL}_n^\vee}$ for $\Lambda$ being the ring of Witt
  vectors of $\overline{\mathbb{F}}_\ell$ with $\ell\neq p$. In our current context ($G$ a connected reductive group
  over $\mathbb{F}_q$), when $G^\ast$ has simply-connected derived subgroup, it is known that
  $\mathsf{B}_{G^\vee}=\mathsf{K}_{G^\ast}$ (see \cite[Thm.\;3.13]{Endomorphism}), so that our main theorem yields the
  equalities
  \[\Lambda\mathsf{E}_G = \Lambda\mathsf{K}_{G^\ast} = \Lambda\mathsf{B}_{G^\vee}\]
  for $\Lambda = \overline{\mathbb{Z}}[\frac{1}{pM}]$. In particular, for $\mathrm{GL}_n$, $M = 1$ and so we provide an
  alternative proof of Helm--Moss's equality.

  {\bf On the proof of the main theorem.}  Identify $ \Lambda\mathsf{E}_G=e_\psi\Lambda G^F e_\psi\subset\Lambda G^F $
  where $e_\psi:=\frac{1}{|U_0^F|}\sum_{u\in U_0^F}\psi(u^{-1})u$ is the primitive central idempotent of $\Lambda U_0^F$
  associated to $\psi$.  We may then consider the {\it symmetrizing form}
  \[
    \tau=\tau_G:=|U_0^F|\mathrm{ev}_{1_{G^F}}:\Lambda\mathsf{E}_G\longrightarrow\Lambda
  \]
  and denote its $\overline{\mathbb{Q}}$-linear extension again by $\tau$. Here $\mathrm{ev}_{1_{G^F}}$ denotes the evaluation map at $1_{G^F}$; recall that a symmetrizing form on a finite
  projective $\Lambda$-algebra $A$ is a map $\tau : A \to \Lambda$ such that the map $(a,b) \mapsto \tau(ab)$ is a
  perfect symmetric bilinear form. It has been
  shown in \cite[Prop.\;2.2]{Endomorphism} that $\tau(\mathsf{K}_{G^\ast})\subset\mathbb{Z}$ and that
  $\tau|_{\Lambda\mathsf{K}_{G^\ast}}:\Lambda\mathsf{K}_{G^\ast}\longrightarrow\Lambda$ is a symmetrizing
  form. Therefore, the equality $\Lambda\mathsf{E}_G=\Lambda\mathsf{K}_{G^\ast}$ will hold if
  \begin{equation}\label{key-lemma}
    \tau(h\pi)\in\Lambda\mbox{\;\;for all\;\;}h\in\Lambda\mathsf{E}_G\mbox{ and }\pi\in\Lambda\mathsf{K}_{G^\ast}.
  \end{equation}
  Indeed, (\ref{key-lemma}) shows that each of $\Lambda\mathsf{E}_G$ and $\Lambda\mathsf{K}_{G^\ast}$ is contained in
  the dual of the other with respect to the above bilinear form; as each is self-dual, they are equal.

After preparations on Deligne-Lusztig characters and Curtis homomorphisms (Section \ref{section-DL-curtis}), we will
  reduce (\ref{key-lemma}) to the study of the condition ``$\tau(h\pi)\in\Lambda$'' for $\pi$ the restriction to $G^{\ast F^{\ast}}$ of a (virtual) algebraic $\overline{\mathbb{F}}_q$-representation of $G^\ast$, by fitting $G^\ast$ into a
  central extension (Section\;\ref{section-DL}) and studying related compatibility questions
  (Sections\;\ref{section-GH1} and \ref{section-GH2}). To study the condition ``$\tau(h\pi)\in\Lambda$'' for such $\pi$, we will
  extend the definition of $\tau(h\pi)$ to $h\in G^F$ (Section\;\ref{section-tau-ext}), reduce the discussion to
  the case where the semisimple part $s$ of $h$ is central in $G$ (Section\;\ref{red-Z-section}), and finally deal with
  the case of central $s$ (Section\;\ref{section-central}).

  {\bf Acknowledgements.} The first author thanks Professor Jean-Fran{\c c}ois Dat, his PhD thesis advisor, for his
  constant support and enlightening opinions on this work. The second author thanks Robert Kurinczuk for bringing the
  work \cite{Endomorphism} of the first author to his attention. We thank Jay Taylor for providing a helpful reference.
  \\[2mm]
  \begin{minipage}{0.10\linewidth}
    \includegraphics[scale=0.06]{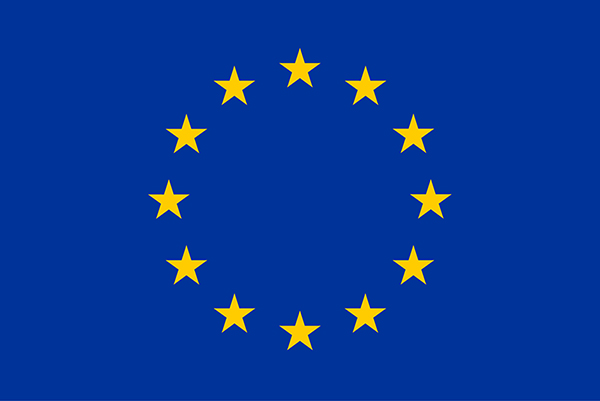}
  \end{minipage}
  \begin{minipage}{0.90\linewidth} {\footnotesize\it This project has received funding from the European Union’s Horizon
      2020 research and \\[-1mm] innovation programme under the Marie Skłodowska-Curie grant agreement No 754362.}
  \end{minipage}

\end{context}

\begin{context}\label{section-DL-curtis} {\bf Preliminaries.}\;--- In this section, we recall some properties of
  Deligne--Lusztig characters and Curtis homomorphisms that we will need later on.

  {\bf Deligne--Lusztig characters.} Let $S$ be an $F$-stable maximal torus of $G$, let $P$ be a Borel subgroup
  containing $S$, and let $V$ be the unipotent radical of $P$. Then we have the Deligne--Lusztig variety (see
  \cite[Def.\;9.1.1]{Digne-Michel})
  \[
    DL_{S \subset P}^G=\{gV\in G/V:g^{-1}F(g)\in V\cdot F(V)\},
  \]
  which admits a (left) $G^F\times (S^F)^{\mathrm{op}}$-action. When there is no need to specify the chosen Borel
  subgroup $P$, we will write $DL_{S \subset P}^G$ simply as $DL_{S}^G$.

  We consider the virtual $\ell$-adic cohomology $H^*_c(\cdot) = \sum_{j \ge 0} (-1)^jH^j_c(\cdot\,, \overline{\mathbb{Q}}_{\ell})$, for $\ell$ a prime
  distinct from $p$. For every character $\chi : S^F \longrightarrow \overline{\mathbb{Q}}^\times$, upon choosing a
  field embedding $\overline{\mathbb{Q}} \hookrightarrow \overline{\mathbb{Q}}_\ell$, we have the corresponding
  Deligne--Lusztig character
  \[
    R_S^G(\chi)(-) := \mathrm{Tr}(-|H_c^\ast(DL_{S\subset P}^G)\otimes_{\overline{\mathbb{Q}}_\ell S^F} \chi) =
    \frac{1}{|S^F|}\sum_{s\in S^F}\mathrm{Tr}((-,s)|H_c^\ast(DL_{S\subset P}^G))\chi(s^{-1}),
  \]
  which is independent of the choice of $P$ and which takes values in $\overline{\mathbb{Q}}_\ell$ {\it {\`a} priori};
  but by \cite[Prop.\;3.3]{Deligne--Lusztig}, for any $(g,s) \in G^F \times S^F$, the trace
  $\mathrm{Tr}((g,s) | H^*_c(DL_{S \subset P}^G))$ is an integer independent of $\ell$, so in fact $R_S^G(\chi)$ takes
  values in $\overline{\mathbb{Q}}$, and it can be verified that ${R}_S^G(\chi)$ is independent of the choices of $\ell$
  and of the embedding $\overline{\mathbb{Q}} \hookrightarrow \overline{\mathbb{Q}}_\ell$.

  {\bf Curtis homomorphisms.} For an $F$-stable maximal torus $S$ of $G$, we consider the Curtis homomorphism
  \[\mathrm{Cur}^G_S : \overline{\mathbb{Q}}\mathsf{E}_G \longrightarrow \overline{\mathbb{Q}}S^F\]
  defined as in \cite[Sec.\;1.7]{Endomorphism} (see also \cite[Thm.\;4.2]{Curtis}). In terms of the Deligne-Lusztig
  dual, the map $\mathrm{Cur}^G_S$ is simply a ``restriction map to a dual torus": indeed, upon fixing an
  $F^\ast$-stable maximal torus $S^\ast$ of $G^\ast$ dual to $S$ (whence a duality
  $\mathrm{Irr}_{\overline{\mathbb{Q}}}(S^F)\simeq S^{\ast F^\ast}$ and thus a ring isomorphism
  $\overline{\mathbb{Q}}S^F\simeq \overline{\mathbb{Q}}^{S^{\ast F^\ast}}$), the map $\mathrm{Cur}^G_S$ is the unique
  ring homomorphism making the following diagram commutative (see \cite[Lem.\;1.6]{Endomorphism}):
  \begin{equation}\label{curtis-torus-restriction}
    \begin{tikzcd}
      \overline{\mathbb{Q}}\mathsf{E}_{G} \arrow[r,rightarrow,"(\ref{EK-identify})"',"\sim"]\arrow[d,rightarrow,"\mathrm{Cur}^G_S"']& \overline{\mathbb{Q}}^{G^{\ast F^\ast}_{\mathrm{ss}}\!/\sim}\arrow[d,rightarrow,"\mathrm{Res}"]\\
      \overline{\mathbb{Q}}S^F \arrow[r,rightarrow,"\sim"]& \overline{\mathbb{Q}}^{S^{\ast F^\ast}}
    \end{tikzcd}
  \end{equation}

  We will later need the following formula of Bonnaf{\'e}--Kessar (\cite[Prop.\;2.5]{Bonnafe-Kessar}, with the missing
  sign factor corrected). For all $h\in\overline{\mathbb{Q}}\mathsf{E}_G\subset\overline{\mathbb{Q}}G^F$,
  \begin{equation}\label{BK-Curtis}
    \mathrm{Cur}_{S}^G(h)=\frac{\epsilon_G\epsilon_S}{|S^F|}\sum_{s\in S^F}\mathrm{Tr}((h,s)|H_c^\ast(DL_{S\subset P}^G))s^{-1}\in\overline{\mathbb{Q}}S^F.
  \end{equation}
  Here, as usual, $\epsilon_G = (-1)^{\mathrm{rk}_{\mathbb{F}_q}(G)}$ for $G$ any reductive group over
  $\mathbb{F}_q$. Observe that (\ref{BK-Curtis}) shows that $\mathrm{Cur}^G_S$ is independent of the choice of $S^\ast$.

\end{context}

\begin{context}\label{section-DL} {\bf On central extensions.}\;--- For our group $G$, we can fit its Deligne--Lusztig
  dual $G^\ast$ into an  $F^\ast$-equivariant exact sequence of reductive groups
  \begin{equation}\label{z-ext}
    1\longrightarrow Z^\ast \longrightarrow H^\ast\longrightarrow G^\ast\longrightarrow 1
  \end{equation}
  where the derived subgroup of $H^\ast$ is simply-connected and $Z^\ast$ is a torus central in $H^\ast$.

  We fix a choice of $F$-equivariant exact sequence of reductive groups
  \begin{equation}\label{z-ext-D}
    1\longrightarrow G\longrightarrow H\xrightarrow{\;\;\kappa\;\;} Z\longrightarrow 1
  \end{equation}
  which is dual to (\ref{z-ext}). Let $T_H$ be an $F$-stable maximal torus of $H$, let $B_H$ be a Borel subgroup of $H$
  containing $T_H$, and let $V$ be the unipotent radical of $B_H$. Then
  \[
    DL_{T_H\subset B_H}^H=\bigsqcup_{z\in Z^F}(DL_{T_H\subset B_H}^H)(z)
  \]
  where for each $z\in Z^F$ we have set
  \[
    (DL_{T_H\subset B_H}^H)(z):=\{hV\in DL_{T_H\subset B_H}^H:\kappa(h)=z\}.
  \]
  Let $T_G=\mathrm{ker}(\kappa|_{T_H}:T_H\twoheadrightarrow Z)$
  (resp.\;$B_G=\mathrm{ker}(\kappa|_{B_H}:B_H\twoheadrightarrow Z)$), which is an $F$-stable maximal torus of $G$
  (resp.\;a Borel subgroup of $G$). Then $T_G\subset B_G$, and the unipotent radical of $B_G$ is also $V$.  As $T_G$ is
  connected, we have $\kappa(T_H^F)=Z^F$, so for each $z\in Z^F$ we may choose a $\dot{z}\in T_H^F$ such that
  $\kappa(\dot{z})=z$. Under the inclusion $G\subset H$, for each $z\in Z^F$ we have
  \[
    DL_{T_H\subset B_H}^H(z)=DL_{T_G\subset B_G}^G \cdot\dot{z}\subset H/V,
  \]
  so that
  \[
    DL_{T_H\subset B_H}^H(z)\simeq DL_{T_G\subset B_G}^G\mbox{ as }(G^F\times (T_G^F)^{\mathrm{op}})\mbox{-varieties}.
  \]

  In terms of virtual $\ell$-adic cohomology we therefore have
  \[
    H_c^\ast(DL_{T_H\subset B_H}^H)=\sum_{z\in Z^F} H_c^\ast(DL_{T_H\subset B_H}^H(z)),
  \]
  and the $H^F\times (T_H^F)^{\mathrm{op}}$-action on $H_c^\ast(DL_{T_H\subset B_H}^H)$ satisfies:
  \[
    \left\lbrace
      \begin{aligned}
        &\mbox{ for every }(h,t)\in H^F\times (T_H^F)^{\mathrm{op}}, (h,t)\cdot H_c^\ast(DL_{T_H\subset B_H}^H(z))\subset H_c^\ast(DL_{T_H\subset B_H}^H(\kappa(ht)z));\\
        &\mbox{ for every }z\in Z^F,\;H_c^\ast(DL_{T_H\subset B_H}^H(z))\simeq H_c^\ast(DL_{T_G\subset B_G}^G)\mbox{ as
        }G^F\times (T_G^F)^{\mathrm{op}}\mbox{-modules}.
      \end{aligned}
    \right.
  \]
  In particular, we obtain the following trace formulae: for $(h,t)\in H^F\times (T_H^F)^{\mathrm{op}}$,
  \begin{equation}\label{tr-formula}
    \left\lbrace
      \begin{aligned}
        \kappa(ht)\neq 1&\Longrightarrow\mathrm{Tr}((h,t)|H_c^\ast(DL_{T_H\subset B_H}^H))=0;\\
        (h,t)\in G^F\times(T_G^F)^{\mathrm{op}}&\Longrightarrow\mathrm{Tr}((h,t)|H_c^\ast(DL_{T_H\subset
          B_H}^H))=|Z^F|\cdot\mathrm{Tr}((h,t)|H_c^\ast(DL_{T_G\subset B_G}^G)).
      \end{aligned}
    \right.
  \end{equation}

  We will later need the compatibility (for $\chi : T_H^F \longrightarrow \overline{\mathbb{Q}}^\times$):
  \begin{equation}\label{tr-formula2}
    R_{T_H}^H(\chi)|_{G^F}=R_{T_G}^G(\chi|_{T_G^F}).
  \end{equation}
  This follows immediately from the defining formula of $R_{T_H}^H(\chi)$ and (\ref{tr-formula}). (See also
  \cite[Prop.\;11.3.10]{Digne-Michel}).
\end{context}

\begin{context}\label{section-GH1} {\bf A compatibility lemma.}\;--- Notation as in Section \ref{section-DL}. We extend
  the $F$-stable Borel subgroup $B_0$ of $G$ in Section \ref{section-intro} (used to determine the Gelfand--Graev module
  $\Gamma_{G,\psi}$) to the $F$-stable Borel subgroup $B_0'$ of $H$, so that $B_0'/B_0=Z$ under (\ref{z-ext-D}); the
  unipotent radical of $B_0'$ is then equal to $U_0$ (the unipotent radical of $B_0$), and the inclusion
  $G^F\subset H^F$ induced by (\ref{z-ext-D}) gives rise to a $\Lambda$-algebra inclusion
  \begin{equation}\label{i-inclusion}
    \Lambda\mathsf{E}_G=e_\psi\Lambda G^F e_\psi\hookrightarrow e_\psi\Lambda H^Fe_\psi = \Lambda\mathsf{E}_{H}.
  \end{equation}
  On the other hand, (\ref{z-ext}) yields the identification
  \begin{equation}\label{ss-ZHG}
    (G^{\ast F^\ast}_{\mathrm{ss}}\!/\!\sim) = (H^{\ast F^\ast}_{\mathrm{ss}}\!/\!\sim)/Z^{\ast F^\ast},
  \end{equation} 
  which enables us to regard functions on $G^{\ast F^\ast}_{\mathrm{ss}}/\!\sim$ as functions on
  $H^{\ast F^\ast}_{\mathrm{ss}}/\!\sim$ which are constant on each $Z^{\ast F^\ast}$-orbit.

  Let us prove the following ``compatibility lemma":

  {\bf Lemma.} {\it The following diagram of rings is commutative:}
  \begin{equation}\label{lemma-diagram}
    \begin{tikzcd}
      \overline{\mathbb{Q}}\mathsf{E}_{G} \arrow[r,rightarrow,"(\ref{EK-identify})"',"\sim"]\arrow[d,hookrightarrow,"(\ref{i-inclusion})"']& \overline{\mathbb{Q}}^{G^{\ast F^\ast}_{\mathrm{ss}}\!/\sim}\arrow[d,hookrightarrow,"(\ref{ss-ZHG})"]\\
      \overline{\mathbb{Q}}\mathsf{E}_{H} \arrow[r,rightarrow,"(\ref{EK-identify})"',"\sim"]&
      \overline{\mathbb{Q}}^{H^{\ast F^\ast}_{\mathrm{ss}}\!/\sim}
    \end{tikzcd}
  \end{equation}

  {\it Proof.} Let $T_G$ and $T_H$ be as in Section \ref{section-DL}, and choose an $F^\ast$-stable maximal torus
  $T_G^\ast$ of $G^\ast$ dual to $T_G$ (resp.\;$T_H^\ast$ of $H$ dual to $T_H$) such that
  ${T_H^\ast/Z^\ast=T_G^\ast}$. Then the Weyl groups of $(G,T_G)$, $(G^\ast,T_G^\ast)$, $(H,T_H)$ and
  $(H^\ast,T_H^\ast)$ are all the same, and we denote this common Weyl group by $W$. For each $w\in W$, choose an
  $F$-stable maximal torus $T_{G,w}$ of $G$ whose $G^F$-conjugacy class corresponds to the $F$-conjugacy class of $w$ in
  $W$ (with respect to $T_G$, so that we may choose $T_{G,1}=T_G$); choose $T_{G,w}^\ast\subset G^\ast$,
  $T_{H,w}\subset H$ and $T_{H,w}^\ast\subset H^\ast$ in a similar way.

  In the toric case where $(G,H)=(T_G,T_H)$, the commutativity of (\ref{lemma-diagram}) follows from toric dualities.

  For the general case of $(G,H)$, we use the Curtis embeddings $\mathrm{Cur}^G=(\mathrm{Cur}_{T_{G,w}}^G)_{w\in W}$ and
  $\mathrm{Cur}^H=(\mathrm{Cur}_{T_{H,w}}^H)_{w\in W}$ (see Section~\ref{section-DL-curtis}) to embed
  (\ref{lemma-diagram}) into the following cubic diagram of rings:
  \begin{equation}\label{lemma-cube}
    \begin{tikzcd}[column sep = small, row sep=small]
      & & \prod\limits_{w\in W} \overline{\mathbb{Q}} T_{G,w}^F\arrow[rr,rightarrow,"\sim"]\arrow[dd,hookrightarrow]&  & \prod\limits_{w\in W}\overline{\mathbb{Q}}^{T_{G,w}^{\ast F^\ast}}\arrow[dd,hookrightarrow]\\
      \overline{\mathbb{Q}}\mathsf{E}_G \arrow[rrr,rightarrow, "\sim" near end, crossing over]\arrow[rru,hookrightarrow,"\mathrm{Cur}^G" near end] \arrow[dd,hookrightarrow]&  &  & \overline{\mathbb{Q}}^{G_{\mathrm{ss}}^{\ast F^\ast}\!/\sim} \arrow[ru,hookrightarrow,"\mathrm{Res}"]&\\
      & & \prod\limits_{w\in W} \overline{\mathbb{Q}} T_{H,w}^F \arrow[rr,rightarrow,"\sim" near start]&  & \prod\limits_{w\in W}\overline{\mathbb{Q}}^{T_{H,w}^{\ast F^\ast}}\\
      \overline{\mathbb{Q}}\mathsf{E}_H \arrow[rrr,rightarrow,"\sim"] \arrow[rru,hookrightarrow,"\mathrm{Cur}^H" near
      end] & & & \overline{\mathbb{Q}}^{H_{\mathrm{ss}}^{\ast F^\ast}\!/\sim}\arrow[from=uu,hookrightarrow,crossing
      over] \arrow[ru,hookrightarrow,"\mathrm{Res}"]&
    \end{tikzcd}
  \end{equation}
  In (\ref{lemma-cube}), the right face is clearly commutative; the top and the bottom faces are commutative by
  (\ref{curtis-torus-restriction}); the back face is the toric case of (\ref{lemma-diagram}) and is hence
  commutative. So to prove the commutativity of (\ref{lemma-diagram}), it remains to show that the left face in
  (\ref{lemma-cube}) is commutative.

  Using (\ref{BK-Curtis}) and the relation $\epsilon_H\epsilon_{T_{H,w}}=\epsilon_G\epsilon_{T_{G,w}}$, the
  commutativity of the left face in (\ref{lemma-cube}) is equivalent to the property that, for all
  $h\in\overline{\mathbb{Q}}\mathsf{E}_G\subset\overline{\mathbb{Q}}G^F$ and all $w\in W$, we have
  \begin{equation}\label{ET-eqX}
    \frac{1}{|T_{H,w}^F|}\sum_{t\in T_{H,w}^F}\mathrm{Tr}((h,t)|H_c^\ast(DL_{T_{H,w}}^H))t^{-1}=\frac{1}{|T_{G,w}^F|}\sum_{t\in T_{G,w}^F}\mathrm{Tr}((h,t)|H_c^\ast(DL_{T_{G,w}}^G))t^{-1}.
  \end{equation}
  By (\ref{tr-formula}) and the fact that $T_{H,w}^F/T_{G,w}^F=Z^F$, we see that (\ref{ET-eqX}) is true for all
  $h\in G^F$, so the left face in (\ref{lemma-cube}) commutes. This completes the proof of the lemma.\qed

\end{context}

\begin{context}\label{section-GH2} {\bf Reduction to the study of $\tau(h\pi_\lambda)$.}\;--- Notation as in
  Section \ref{section-DL}. As $Z^{\ast F^\ast}$ is central in $H^{\ast F^\ast}$, the association of each irreducible
  $\overline{\mathbb{F}}_qH^{\ast F^\ast}$-module to its restriction to $Z^{\ast F^\ast}$ 
  induces a $\widehat{Z^{\ast F^{\ast}}}$-graded decomposition
  \begin{equation}\label{K-grading}
    \mathsf{K}_{H^\ast}=\bigoplus\limits_{\lambda\in \widehat{Z^{\ast F^{\ast}}}}(\mathsf{K}_{H^\ast})_\lambda\mbox{\;\;with\;\;}\mathsf{K}_{G^\ast}=(\mathsf{K}_{H^\ast})_{1}.
  \end{equation}
  In particular, we have a ring inclusion $\mathsf{K}_{G^{\ast}}\subset \mathsf{K}_{H^\ast}$, and it is evident that the
  following diagram of rings is commutative (where $\mathrm{br}$ denotes the Brauer character map):
  \begin{equation}\label{lemma-diagram2}
    \begin{tikzcd}
      \overline{\mathbb{Q}}\mathsf{K}_{G^\ast}\arrow[d,hookrightarrow,"(\ref{K-grading})"']\arrow[r,rightarrow,"\mathrm{br}","\sim"'] &   \overline{\mathbb{Q}}^{G^{\ast F^\ast}_{\mathrm{ss}}\!/\sim} \arrow[d,hookrightarrow,"(\ref{ss-ZHG})"]\\
      \overline{\mathbb{Q}}\mathsf{K}_{H^\ast} \arrow[r,rightarrow,"\mathrm{br}","\sim"'] &
      \overline{\mathbb{Q}}^{H^{\ast F^\ast}_{\mathrm{ss}}\!/\sim}
    \end{tikzcd}
  \end{equation}

Let $h\in\Lambda\mathsf{E}_G$ and $\pi\in\mathsf{K}_{G^\ast}$. Via the commutative diagrams (\ref{lemma-diagram}) and (\ref{lemma-diagram2}), we can define the product $h\pi$ consistently as an   element of $\overline{\mathbb{Q}}\mathsf{E}_G$, $\overline{\mathbb{Q}}\mathsf{E}_H$,  $\overline{\mathbb{Q}}^{G_{\mathrm{ss}}^{\ast F^\ast}\!/\sim}$,  $\overline{\mathbb{Q}}^{H_{\mathrm{ss}}^{\ast F^\ast}\!/\sim}$, $\overline{\mathbb{Q}}\mathsf{K}_{G^\ast}$ or  $\overline{\mathbb{Q}}\mathsf{K}_{H^\ast}$. 
As 
 $
  \tau_G=|U_0^F|\mathrm{ev}_{1_{G^F}}=|U_0^F|\mathrm{ev}_{1_{H^F}}=\tau_H,
 $ we deduce that
  \begin{equation}\label{pi-lambda-red}
    \tau_G(h\pi)=\tau_H(h\pi).
  \end{equation}
  Therefore, if we can prove (\ref{key-lemma}) for $\tau_H(h\pi)$, then we can prove it for $\tau_G(h\pi)$.

Now let $\mathsf{K}(G^\ast\mbox{-mod})$ be the Grothendieck ring of the category of
  finite-dimensional algebraic $G^\ast$-modules and let $\mathsf{K}^\circ_{G^\ast}$ be the image of the restriction map
  \[\mathrm{Res}:\mathsf{K}(G^\ast\mbox{-mod})\longrightarrow\mathsf{K}_{G^\ast}.\]
  Adopting similar notation for $H$, we have that $\mathrm{Res}$ is surjective (\cite[Thm.\;7.4]{Steinberg} and \cite[Thm.\;3.10]{Herzig}) so that $\mathsf{K}^\circ_{H^\ast} = \mathsf{K}_{H^\ast}$. We are therefore reduced to proving that $\tau_H(h\pi) \in \Lambda$ for $\pi \in \mathsf{K}^\circ_{H^\ast}$.
This turns out to be true without the assumption that $H^\ast$ has simply-connected derived subgroup; in the following, we shall thus return to the group $G$ and study the condition ``$\tau_G(h\pi)\in\Lambda$ for $\pi \in \mathsf{K}^\circ_{G^\ast}$.''  
\end{context}

\begin{context}\label{section-tau-ext} {\bf An extension $\widetilde{\tau}$ for $\tau(h\pi)$.}\;---
  We return to the group $G$ (the derived subgroup of $G^\ast$ may not be simply-connected) and write $\tau_G=\tau$.
  Let $T$ be an $F$-stable maximal torus of $G$, let $W=N_G(T)/T$ be the Weyl group of $(G,T)$, and let $T_w$ be an
  $F$-stable maximal torus of $G$ associated with $w\in W$ (with respect to $T$) as in the proof of
  (\ref{lemma-diagram}).  Recall the identification
  $\overline{\mathbb{Q}}\mathsf{E}_G=\overline{\mathbb{Q}}\mathsf{K}_{G^\ast}$ from (\ref{EK-identify}). Then, for
  $h\in\overline{\mathbb{Q}}\mathsf{E}_G$ and $\pi\in\overline{\mathbb{Q}}\mathsf{K}_{G^\ast}$:
  \begin{align*}
    \tau(h\pi)&=\frac{1}{|W|}\sum_{w\in W}\mathrm{ev}_{1_{T_w^F}}(\mathrm{Cur}_{T_w}^G(h\pi))\quad(\mbox{by \cite[Eq.\;3.5]{Bonnafe-Kessar}}) \\
              &=\frac{1}{|W|}\sum_{w\in W}\mathrm{ev}_{1_{T_w^F}}(\mathrm{Cur}_{T_w}^G(h)\cdot\mathrm{Cur}_{T_w}^G(\pi))  \quad\mbox{($\mathrm{Cur}_{T_w}^G$ is a ring homomorphism)} \\
              &=\frac{1}{|W|}\sum_{w\in W}\sum_{t\in T_w^F}\mathrm{Cur}_{T_w}^G(h)(t^{-1})\cdot\mathrm{Cur}_{T_w}^G(\pi)(t)\\
              &=\frac{1}{|W|}\sum_{w\in W}\sum_{t\in T_w^F}\frac{\epsilon_G\epsilon_{T_w}}{|T_w^F|}\cdot\mathrm{Tr}((h,t)|H_c^\ast(DL_{T_w}^G))\cdot\mathrm{Cur}_{T_w}^G(\pi)(t) \quad\mbox{(by (\ref{BK-Curtis}))} \\
              &=\frac{1}{|W|}\sum_{w\in W}\frac{\epsilon_G\epsilon_{T_w}}{|T_w^F|}\sum_{t\in T_w^F}\sum_{\chi\in\mathrm{Irr}_{\overline{\mathbb{Q}}}(T_w^F)}R_{T_w}^G(\chi)(h)\cdot\chi(t)\cdot\mathrm{Cur}_{T_w}^G(\pi)(t)\quad\mbox{(trace formula)}\\
              &=\frac{1}{|W|}\sum_{w\in W}\frac{\epsilon_G\epsilon_{T_w}}{|T_w^F|}\sum_{\chi\in\mathrm{Irr}_{\overline{\mathbb{Q}}}(T_w^F)}R_{T_w}^G(\chi)(h)\cdot\chi(\mathrm{Cur}_{T_w}^G(\pi)). \stepcounter{equation}\tag{\theequation}\label{BK-Curtis-2a} 
  \end{align*}

  Using the formula (\ref{BK-Curtis-2a}), we can extend the $\overline{\mathbb{Q}}$-bilinear map
  \[
    \overline{\mathbb{Q}}\mathsf{E}_G\times\overline{\mathbb{Q}}\mathsf{K}_{G^\ast}\longrightarrow\overline{\mathbb{Q}},\quad(h,\pi)\longmapsto
    \tau(h\pi),
  \]
  to a $\overline{\mathbb{Q}}$-bilinear map
  $\widetilde{\tau}(\cdot,\cdot):\overline{\mathbb{Q}}G^F\times\overline{\mathbb{Q}}\mathsf{K}_{G^\ast}\longrightarrow\overline{\mathbb{Q}}$
  by setting, for $h\in \overline{\mathbb{Q}}G^F$ and $\pi\in\overline{\mathbb{Q}}\mathsf{K}_{G^\ast}$,
  \begin{equation}\label{BK-Curtis-R}
    \widetilde{\tau}(h,\pi):=\frac{1}{|W|}\sum_{w\in W}\frac{\epsilon_G\epsilon_{T_w}}{|T_w^F|}\sum_{\chi\in\mathrm{Irr}_{\overline{\mathbb{Q}}}(T_w^F)}R_{T_w}^G(\chi)(h)\cdot\chi(\mathrm{Cur}_{T_w}^G(\pi)).
  \end{equation} 
  We then have
  \begin{equation}\label{rest-tau-A}
    \tau(h\pi)=\widetilde{\tau}(h,\pi)\;\mbox{for all}\; h\in\overline{\mathbb{Q}}\mathsf{E}_G\;\mbox{and all}\;\pi\in\overline{\mathbb{Q}}\mathsf{K}_{G^\ast}.
  \end{equation}

  The formula (\ref{BK-Curtis-R}) for $\widetilde{\tau}$ involves choices of $T$ and $T_w$; we now derive an intrinsic
  formula for $\widetilde{\tau}$ as follows.

  Let $\mathcal{T}_G$ be the set of $F$-stable maximal tori of $G$, and let $\mathcal{T}_G/G^F$ be the set of
  $G^F$-conjugacy classes in $\mathcal{T}_G$. For each $S\in\mathcal{T}_G$, let $W_G(S)=N_G(S)/S$. Since the isomorphism
  class of $T_w$ depends only on the $F$-twisted conjugacy class of $w \in W$, and the stabiliser of $w \in W$ under
  $F$-twisted conjugacy may be identified with $W_G(T_w)^F$, we have that there are $\frac{|W|}{|W_G(S)^F|}$ elements
  $w \in W$ such that $T_w$ is $G^F$-conjugate to $S$.  By (\ref{BK-Curtis-R}), for
  $h\in \overline{\mathbb{Q}}G^F$ and $\pi\in\overline{\mathbb{Q}}\mathsf{K}_{G^\ast}$, we have:
  \begin{align*}
    \widetilde{\tau}(h,\pi)&=\sum_{S\in\mathcal{T}_G/G^F}\frac{\epsilon_G\epsilon_S}{|W_G(S)^F|}\frac{1}{|S^F|}\sum_{\chi\in\mathrm{Irr}_{\overline{\mathbb{Q}}}(S^F)}R_S^G(\chi)(h)\cdot\chi(\mathrm{Cur}_S^G(\pi))\\
                           &=\frac{1}{|G^F|}\sum_{S\in\mathcal{T}_G}\epsilon_G\epsilon_S\sum_{\chi\in\mathrm{Irr}_{\overline{\mathbb{Q}}}(S^F)}R_S^G(\chi)(h)\cdot\chi(\mathrm{Cur}_S^G(\pi)).\stepcounter{equation}\tag{\theequation}\label{BK-Curtis-I}
  \end{align*}

\end{context}

\begin{context}\label{red-Z-section} {\bf Reduction to the case of central $s$.}\;--- From now on, let $h=su\in G^F$
  with $s\in G^F$ (resp.\;$u\in G^F$) the semisimple (resp.\;unipotent) part in the Jordan decomposition of $h$. Recall
  Deligne--Lusztig's character formula \cite[Thm.\;4.2]{Deligne--Lusztig} for each $F$-stable maximal torus $S$ of $G$:
  (notation: $\mathrm{ad}(g)x={}^gx=gxg^{-1}$)
  \begin{equation}\label{DL-character} R_{S}^G(\chi)(h)=\frac{1}{|C_G(s)^{\circ F}|}\sum_{\substack{g\in G^F\\
        g^{-1}sg\in
        S^F}}Q_{\mathrm{ad}(g)S}^{C_G(s)^\circ}(u)\cdot({}^g\chi)({s})
  \end{equation} where  $Q_S^G=R_S^G(\mathbf{1})|_{G_{\mathrm{unip}}^F}$ denotes the Green function and 
  $C_G(s)^\circ$ is the identity component of the centralizer of $s$ in $G$.

  We shall write $\widetilde{\tau}=\widetilde{\tau}_G$ to specify the group $G$. Substituting (\ref{DL-character}) into
  (\ref{BK-Curtis-I}), we obtain: (below, $\pi\in\mathsf{K}_{G^\ast}$)
  \begin{align*}
    \widetilde{\tau}_G(h,\pi)&=\frac{1}{|G^F|}\sum_{S\in\mathcal{T}_G}\epsilon_G\epsilon_S\sum_{\chi\in\mathrm{Irr}_{\overline{\mathbb{Q}}}(S^F)}\frac{1}{|C_G(s)^{\circ
                               F}|}\sum_{\substack{g\in G^F\\ g^{-1}sg\in
    S^F}}Q_{\mathrm{ad}(g)S}^{C_G(s)^\circ}(u)\cdot\chi({}^{g^{-1}}s\cdot\mathrm{Cur}_S^G(\pi))\\
                             &=\frac{1}{|G^F|}\frac{1}{|C_G(s)^{\circ F}|}\sum_{S\in\mathcal{T}_G}\epsilon_G\epsilon_S|S^F|\sum_{\substack{g\in G^F\\
    g^{-1}sg\in S^F}}Q_{\mathrm{ad}(g)S}^{C_G(s)^\circ}(u)\cdot\mathrm{Cur}_{S}^G(\pi)({}^{g^{-1}}(s^{-1}))\\
  \intertext{(where we have applied the orthogonality of characters)}
   &=\frac{1}{|G^F|}\frac{1}{|C_G(s)^{\circ F}|}\sum_{g\in G^F}\sum_{\substack{S\in\mathcal{T}_G\\
    s\in(\mathrm{ad}(g)S)^F}}\epsilon_G\epsilon_S|S^F|\cdot
    Q_{\mathrm{ad}(g)S}^{C_G(s)^\circ}(u)\cdot\mathrm{Cur}_{\mathrm{ad}(g)S}^G(\pi)(s^{-1})\\
  \intertext{(where we have used $\mathrm{Cur}_{\mathrm{ad}(g)S}^G(\pi)({}^gx) = \mathrm{Cur}_S^G(\pi)(x)$ for $g\in
  G^F$)}
  &=\frac{1}{|C_G(s)^{\circ F}|}\sum_{\substack{S\in\mathcal{T}_G\\ s\in S^F}}\epsilon_G\epsilon_S|S^F|\cdot Q_{S}^{C_G(s)^\circ}(u)\cdot\mathrm{Cur}_{S}^G(\pi)(s^{-1})\quad(S\longmapsto
    \mathrm{ad}(g^{-1})S)\\ 
    &=\frac{1}{|C_G(s)^{\circ F}|}\sum_{\substack{S\in\mathcal{T}_{C_G(s)^\circ}\\s\in S^F}}\epsilon_G\epsilon_S|S^F|\cdot
    Q_{S}^{C_G(s)^\circ}(u)\cdot\mathrm{Cur}_{S}^G(\pi)(s^{-1}),\stepcounter{equation}\tag{\theequation}\label{tilde-tau-I}
  \end{align*} where the last equality holds because for $S\in\mathcal{T}_G$, if $S^F$ contains $s$ then $S\subset
  C_G(s)^\circ$.

Recall the subring $\mathsf{K}^\circ_{G^\ast} \subset \mathsf{K}_{G^\ast}$ from section~\ref{section-GH2}.

{\bf Lemma.} {\it Let $\Lambda_0$ be a subring of  
$\overline{\mathbb{Q}}$. Fix an $h=su\in G^F$ as above, and consider the following statement: 
  \begin{equation}\label{weak-ver} 
  \widetilde{\tau}_G(h,\pi)\in\Lambda_0\quad\mbox{ for all }\pi\in\mathsf{K}_{G^\ast}^\circ.
  \end{equation}
Suppose that (\ref{weak-ver}) is true when $G$ therein is replaced by $C_G(s)^\circ$ (by \cite[Prop.\;3.5.3]{Digne-Michel}, $u\in C_G(s)^\circ$ and hence  $h\in C_G(s)^{\circ F}$). Then (\ref{weak-ver}) is true for $G$.}
  
As $s$ is central in $C_G(s)^\circ$, 
this lemma will reduce the study of the condition (\ref{weak-ver}) to the case where $s$ is central in $G$.

  {\it Proof of lemma.} 
  First, we require a certain special set of generators for $\mathsf{K}^\circ_{G^\ast}$. As shown in \cite[Ch.\;II.2]{Jantzen}, for every maximal torus $T^\ast$ of $G^\ast$, the associated
  formal character map gives a ring isomorphism
  \[
    \mathrm{ch}:\mathsf{K}(G^\ast\mbox{-mod})\xrightarrow{\;\;\sim\;\;}\mathbb{Z}[X(T^\ast)]^W
  \]
  where $X(T^\ast)=\mathrm{Hom}_{\mathrm{alg}}(T^\ast,\mathbb{G}_m)$ is the character group of $T^\ast$ and $W$ is the
  Weyl group of $(G^\ast,T^\ast)$. For $\lambda\in X(T^\ast)$, set 
  \[
    r_{G,\lambda}:=\sum\limits_{\mu\in W\lambda}\mu\in\mathbb{Z}[X(T^\ast)]^W\mbox{\quad
      and\quad}\pi_{G,\lambda}:=\mathrm{ch}^{-1}(r_{G,\lambda})|_{G^{\ast F^\ast}}\in\mathsf{K}^\circ_{G^\ast}, 
  \]
  where for $\lambda\in X(T^\ast)$, $W\lambda$ denotes the $W$-orbit of $\lambda$. Note that the $\mathbb{Z}$-module
  $\mathbb{Z}[X(T^\ast)]^W$ is generated by $\{r_{G,\lambda}:\lambda\in X(T^\ast)\}$, and so the $\pi_{G,\lambda}$ generate $\mathsf{K}^\circ_{G^\ast}$ as a $\mathbb{Z}$-module.

Choose an $F$-stable maximal torus $T$ of $G$ containing $s$, 
so that $T$ is  also an $F$-stable maximal torus of $C_G(s)^\circ$.
To verify (\ref{weak-ver}) for the chosen $h$, it suffices to show that  $\widetilde{\tau}_G(h,\pi_{G,\lambda})\in\Lambda_0$ for all $\lambda\in X(T^\ast)$.

  Let $S$ be an $F$-stable maximal torus of $G$, 
  choose an $F^\ast$-stable maximal torus
  $S^\ast$ of $G^\ast$ dual to $S$ and with a duality
  $\widehat{\cdot}:S^F\xrightarrow{\;\sim\;}\mathrm{Irr}_{\overline{\mathbb{F}}_q}(S^{\ast F^\ast})$, and fix a choice
  of $g\in G^\ast$ such that $S^\ast={}^gT^\ast$. This duality and the fixed embedding
  $\overline{\mathbb{F}}^\times_q \hookrightarrow \overline{\mathbb{Q}}^\times$ allow us to identify
  $\overline{\mathbb{Q}}S^F$ with $\overline{\mathbb{Q}}^{S^{\ast F^\ast}}$.  For each $\mu\in X(T^\ast)$, set
  $\mu_{S^\ast}={}^g\mu\in X(S^\ast)$, and define $\phi_S(\mu)\in S^F$ by the relation
  $\mu_{S^\ast}|_{S^{\ast F^\ast}}=\widehat{\phi_S(\mu)}\in\mathrm{Irr}_{\overline{\mathbb{F}}_q}(S^{\ast F^\ast})$. We
  then have a map $ \phi_S:X(T^\ast)\longrightarrow S^F $ which extends to a ring homomorphism
  \[ \phi_S:\overline{\mathbb{Q}}[X(T^\ast)]\longrightarrow \overline{\mathbb{Q}}S^F = \overline{\mathbb{Q}}^{S^{\ast
        F^\ast}}.
  \] The following diagram then commutes (where $W=W_{G^\ast}(T^\ast)=N_{G^\ast}(T^\ast)/T^\ast$): 
  \[
    \begin{tikzcd}
      \overline{\mathbb{Q}}\mathsf{K}(G^\ast\mbox{-mod})\arrow[rr,rightarrow,"\mathrm{ch}","\sim"']\arrow[d,"\mathrm{Res}_{G^{\ast
          F^\ast}}^{G^\ast}"'] & &\overline{\mathbb{Q}}[X(T^\ast)]^W \arrow[d,rightarrow,"\phi_S"]\\
      \overline{\mathbb{Q}}\mathsf{K}_{G^\ast} \arrow[r,rightarrow,"\mathrm{br}","\sim"'] &
      \overline{\mathbb{Q}}^{G^{\ast F^\ast}} \arrow[r,rightarrow,"\mathrm{Res}^{G^{\ast F^{\ast}}}_{S^{\ast F^\ast}}"]
      & \overline{\mathbb{Q}}^{S^{\ast F^\ast}}
    \end{tikzcd}
  \] Combining this with (\ref{curtis-torus-restriction}) we see that the following diagram of rings also commutes:
\begin{equation}\label{phi-compatible}
    \begin{tikzcd} \overline{\mathbb{Q}}\mathsf{K}(G^\ast\mbox{-mod})\arrow[r,"\mathrm{Res}_{G^{\ast
          F^\ast}}^{G^\ast}"]\arrow[d,rightarrow,"\mathrm{ch}"',"\sim" sloped]&
      \overline{\mathbb{Q}}\mathsf{K}_{G^\ast}\arrow[r,equal,"(\ref{EK-identify})"]&
      \overline{\mathbb{Q}}\mathsf{E}_G\arrow[r,"\mathrm{Cur}_S^G"]& \overline{\mathbb{Q}}S^F\\
      \overline{\mathbb{Q}}[X(T^\ast)]^{W}\arrow[rrru,"\phi_S"']
    \end{tikzcd}
  \end{equation}

  The commutative diagram (\ref{phi-compatible}) gives the relation
  \begin{equation}\label{phi-keyprop} \mathrm{Cur}_{S}^G(\pi_{G,\lambda})=\phi_S(r_{G,\lambda}).
  \end{equation} 
    Via the identifications
  \[ W_{C_G(s)^{\circ \ast}}(T^\ast)=W_{C_G(s)^\circ}(T)\leq W_G(T)=W_{G^\ast}(T^\ast),
  \] we may write
  $W_{G^\ast}(T^\ast)\lambda=\bigsqcup\limits_{\lambda'\in \Omega} W_{C_G (s)^{\circ\ast}}(T^\ast)\lambda'$ for some
  finite subset $\Omega$ of $W_{G^\ast}(T^\ast)\lambda$, so that
  $r_{G,\lambda}=\sum\limits_{\lambda'\in\Omega}r_{C_G(s)^\circ,\lambda'}$ and then
  $\mathrm{Cur}_S^G(\pi_{G,\lambda})=\sum\limits_{\lambda'\in\Omega}\mathrm{Cur}_{S}^{C_G(s)^\circ}(\pi_{C_G(s)^\circ,\lambda'})$
  by (\ref{phi-keyprop}). Applying (\ref{tilde-tau-I}) to $\pi=\pi_{G,\lambda}$, we thus deduce that
  \begin{equation}\label{tilde-tau-3}
    \widetilde{\tau}_G(h,\pi_{G,\lambda})=\epsilon_G\epsilon_{C_G(s)^\circ}\sum_{\lambda'\in\Omega}\widetilde{\tau}_{C_G(s)^\circ}(h,\pi_{C_G(s)^\circ,\lambda'}).
  \end{equation}
  By (\ref{tilde-tau-3}) and the assumption of the lemma, 
  we get $\widetilde{\tau}_G(h,\pi_{G,\lambda})\in\Lambda_0$ for all $\lambda\in X(T^\ast)$, whence
  $\widetilde{\tau}_G(h,\pi)\in\Lambda_0$ for all $\pi\in\mathsf{K}_{G^\ast}^\circ$. \qed

\end{context}

\begin{context}\label{section-central} {\bf The case of central $s$.}\;--- Keep the notation $T$, $W$ and $T_w$ as in
  Section \ref{section-tau-ext}. Let $\pi \in \mathsf{K}_{G^\ast}$, let $h=su\in G^F$ be as in Section \ref{red-Z-section}, and suppose furthermore that
  $s$ lies in the centre of $G$. Then $C_G(s)^\circ=G$, and (\ref{DL-character}) becomes
  $R_S^G(\chi)(h)=Q_S^G(u)\chi(s)$, so that (\ref{BK-Curtis-R}) yields
  \begin{align*} \widetilde{\tau}_G(h,\pi)&=\frac{1}{|W|}\sum_{w\in
                                                    W}\frac{\epsilon_G\epsilon_{T_w}}{|T_w^F|}\sum_{\chi\in\mathrm{Irr}_{\overline{\mathbb{Q}}}(T_w^F)}Q_{T_w}^G(u)\cdot\chi(s)\cdot\chi(\mathrm{Cur}_{T_w}^G(\pi))\\
                                                  &=\frac{1}{|W|}\sum_{w\in W}\epsilon_G\epsilon_{T_w}
                                                    Q_{T_w}^G(u)\langle{\pi}|_{T_w^{\ast
                                                    F^\ast}},\widehat{s^{-1}}\rangle_{T_w^{\ast F^\ast}}\quad(\mbox{orthogonality of characters})\\
                                                  &=\langle{\tilde{\pi}},\gamma\rangle_{G^{\ast
                                                    F^\ast}}\stepcounter{equation}\tag{\theequation}\label{central-eqA}
  \end{align*} where (using \cite[Prop.\;7.3]{Deligne--Lusztig})
  \[ \gamma:=\frac{1}{|W|}\sum_{w\in W}\epsilon_G\epsilon_{T_w} Q_{T_w}^G(u)\mathrm{Ind}_{T_w^{\ast F^\ast}}^{G^{\ast
        F^\ast}}\widehat{s^{-1}}=\frac{1}{|W|}\sum_{w\in
      W}Q_{T_w}^G(u)R_{T_w^\ast}^{G^\ast}\widehat{s^{-1}}\otimes\mathrm{St}_{G^\ast}
  \] 
and $\tilde{\pi}$ is any extension of the Brauer character $\pi$ to an ordinary virtual character (which exists by \cite[Thm.\;33]{Serre}). 
  As $s$ lies in the centre of $G$, $\widehat{s^{-1}}$ is in fact a multiplicative character of $G^{\ast F^\ast}$, so
  \begin{align} \gamma &=\gamma'\otimes\widehat{s^{-1}}\otimes\mathrm{St}_{G^\ast} \nonumber \\ \intertext{with} \gamma'
    &:=\frac{1}{|W|}\sum_{w\in W}Q_{T_w}^G(u)R_{T_w^\ast}^{G^\ast}(\mathbf{1}).\label{central-eqB}
  \end{align}
Our strategy will be to show that $\gamma'$ is a $\overline{\mathbb{Q}}$-linear combination of irreducible $G^{\ast F^\ast}$-representations with only bad primes appearing in the denominators.

We need some facts from the theory of almost characters, following \cite[Ch.\;3-4]{Lusztig}; in the notation of that
  book, we are considering the case $n = 1$ and $L$ trivial. See also \cite[Sec.\;7.3]{Carter} for a concise exposition,
  but with some extraneous hypotheses. Let $c$ be the order of the automorphism $F$ on $W$ (when $G$ is split, we have
  $c=1$); denote by $\widehat{W}_{\mathrm{ex}}$ the set of all $\phi\in\mathrm{Irr}_{\mathbb{Q}}(W)$ which can be
  extended to a $\mathbb{Q}$-valued irreducible character of $\widehat{W}:=W\rtimes\langle F\rangle$ (by
  \cite[Cor.\;1.15]{Springer}, every irreducible representation of $W$ over a characteristic 0 field is defined over
  $\mathbb{Q}$); for each $\phi\in\widehat{W}_{\mathrm{ex}}$, there exists such an extension (in fact, exactly two)
  $\widetilde{\phi}\in\mathrm{Irr}_{\mathbb{Q}}(\widehat{W})$. Fixing a choice of such $\widetilde{\phi}$, we then call 
  \[ R_{\widetilde{\phi}}^{G^\ast}:=\frac{1}{|W|}\sum_{w\in W}\widetilde{\phi}(wF)R_{T_w^\ast}^{G^\ast}(\mathbf{1})
  \] an {\it almost character} of $G^{\ast F^\ast}$.

  Recall from Section\;\ref{section-intro} the definition of bad primes for $G$. Note that a prime is bad for $G$ if and
  only if it is bad for $G^*$. Define
  \begin{equation}\label{def-M-G} \mbox{$M_G=$ product of all bad primes for $G$}.
  \end{equation} Using Lusztig's work on unipotent characters,
  \begin{equation}\label{almost-eq1}
    \begin{aligned} &\mbox{each almost character $R_{\widetilde{\phi}}^{G^\ast}$ is a
        $\overline{\mathbb{Z}}[\frac{1}{M_G}]$-linear combination} \\[-1mm] &\mbox{of irreducible $\overline{\mathbb{Q}}$-valued
        unipotent characters of $G^{\ast F^\ast}$}.
    \end{aligned}
  \end{equation} Indeed, if $G^\ast$ has connected centre, then \cite[Thm.\;4.23]{Lusztig} expresses
  $R_{\widetilde{\phi}}$ as a linear combination of unipotent characters of $G^{\ast F^\ast}$. By
  \cite[(4.21.7)]{Lusztig}, the denominators divide the orders of certain groups $\mathcal{G}_{\mathcal{F}}$ of the form
  $\prod\mathcal{G}_{\mathcal{F}_i}$ where the product is over the irreducible factors of the root system of
  $G^\ast$. Each $\mathcal{G}_{\mathcal{F}_i}$ is defined in a case-by-case fashion, in a way depending only on the
  corresponding irreducible factor of the root system, in \cite[4.4--4.13]{Lusztig}, and has order divisible only by bad
  primes for that factor. If $G^\ast$ does not have connected centre then we choose a short exact sequence
  \[1 \to G^\ast \to H^\ast \to Z^\ast \to 1\] as in (\ref{z-ext-D}) (with the roles of $G^\ast$ and $G$
  reversed). Extending the chosen maximal $F^\ast$-stable torus and Borel from $G^\ast$ to $H^\ast$ as in Section
  \ref{section-DL}, we may identify the Weyl groups of $G^\ast$ and $H^\ast$. Using (\ref{tr-formula2}) (with
  $\chi=\mathbf{1}$ therein), we then have
  \[R_{\widetilde{\phi}}^{H^\ast}|_{G^{\ast F^\ast}} = R_{\widetilde{\phi}}^{G^\ast},\] whence
  $R_{\widetilde{\phi}}^{G^\ast}$ is $\overline{\mathbb{Z}}[\frac{1}{M_G}]$-linear combination of restrictions to
  $G^{\ast F^\ast}$ of unipotent characters of $H^{\ast F^\ast}$. However, the restriction to $G^{\ast F^\ast}$ of a
  unipotent character of $H^{\ast F^\ast}$ is a unipotent character by \cite[Prop.\;11.3.8]{Digne-Michel}, so
  (\ref{almost-eq1}) follows.

  Now we prove the following lemma:

  {\bf Lemma.} {\it The sum  
  $\gamma'$ in (\ref{central-eqB}) is a finite
    $\overline{\mathbb{Z}}[\frac{1}{M_G}]$-linear combination of almost characters of $G^{\ast F^\ast}$.  }

  {\it Proof.} We have  
  \begin{align*} \gamma'&=\frac{1}{|W|}\sum_{w\in W}R_{T_w}^G(\mathbf{1})(u)R_{T_w^\ast}^{G^\ast}(\mathbf{1})\\
                        &=\frac{1}{|W|}\sum_{w\in
                          W}R_{T_w}^G(\mathbf{1})(u)\sum_{\phi\in\widehat{W}_{\mathrm{ex}}}\widetilde{\phi}(wF)R_{\widetilde{\phi}}^{G^\ast}
                          \quad(\mbox{see \cite[p.\;76]{Carter}})\\
                        &=\sum_{\phi\in\widehat{W}_{\mathrm{ex}}}R_{\widetilde{\phi}}^G(u)R_{\widetilde{\phi}}^{G^\ast};
  \end{align*} by (\ref{almost-eq1}) and the fact that character values of representations of finite groups are algebraic
  integers, all $R_{\widetilde{\phi}}^G$ must take values in $\overline{\mathbb{Z}}[\frac{1}{M_G}]$.\qed

  {\bf Remark.} In the above lemma, the class function $\gamma'$ can in fact be written as a finite
  $\overline{\mathbb{Z}}$-linear combination of almost characters of $G^{\ast F^\ast}$. We won't need this stronger property of $\gamma'$ later, so here we only briefly explain how to achieve this, following the complete proof in \cite[Rmk.\;of Lem.\;2.23]{Thesis}. First, one uses 
a theorem of Shoji (\cite[Thm.\;5.5]{Shoji}; see also \cite[Thm.\;13.2.3]{Digne-Michel}) to get that
  $ Q_{T_w}^G(u)=\mathrm{Tr}(wF|H_c^\ast(\mathcal{B}_u))$ for all $w\in W$, where $\mathcal{B}_u$ is the variety of
  Borel subgroups of $G$ containing $u$. One then studies the contribution of each composition factor $V$ of the finite-dimensional $\overline{\mathbb{Q}}_\ell\widehat{W}$-module $H_c^\ast(\mathcal{B}_u)$ ($\ell\neq p$) to $\mathrm{Tr}(wF|H_c^\ast(\mathcal{B}_u))$; one proves that $\mathrm{Tr}(wF|V)\neq 0$ only if $V|_W$ is irreducible, and in this case $\mathrm{Tr}(wF|V)=\chi_V(F)\cdot\mathrm{Tr}(wF|\widetilde{\phi})$ for some linear character $\chi_V:\langle F\rangle\longrightarrow\overline{\mathbb{Q}}_\ell^\times $ and some $\widetilde{\phi}\in\mathrm{Irr}_{\mathbb{Q}}(\widehat{W})$ fitting the definition of the almost character $R_{\widetilde{\phi}}^{G^\ast}$ and on which $F^c$ acts trivially, so that $\gamma'$ is the sum of finitely many $\chi_V(F)\cdot R_{\widetilde{\phi}}^{G^\ast}$ with $V|_W$ irreducible. As all eigenvalues of  the endomorphism $F$ on $H_c^\ast(\mathcal{B}_u)$ lie in $\overline{\mathbb{Z}}$ (see \cite[Lem.\;1.7]{Deligne}), each $\chi_V(F)$ must lie in $\overline{\mathbb{Z}}^\times$, so $\gamma'$ is a finite
  $\overline{\mathbb{Z}}$-linear combination of almost characters of $G^{\ast F^\ast}$, as desired.

  Using the previous lemma, (\ref{central-eqA}), (\ref{central-eqB}) and (\ref{almost-eq1}), we get the following
  proposition:

  {\bf Proposition.} {\it We have $\widetilde{\tau}_G(h,\pi)\in\overline{\mathbb{Z}}[\frac{1}{M_G}]$ for all
    $\pi \in \mathsf{K}_{G^\ast}$ and all $h\in G^F$ whose semisimple part $s$ is central in $G$. ($M_G$ is as in
    (\ref{def-M-G}).)}

  {\bf End of proof of the main theorem in Section \ref{section-intro}.} From now on, we remove the assumption that $s$
  is central in $G$.

  Observe that a prime number that is bad for $C_G(x)^\circ$ with $x$ a semisimple element of $G^F$ is also bad for $G$;
  indeed, this follows from the definition of bad primes in Section \ref{section-intro} and from the following two
  facts: (i) if $G$ is simple of type $A$ (resp.\;of classical type), then the centralizer of every semisimple element
  of $G$ has only factors of type $A$ (resp.\;of classical type); (ii) if $G$ is simple of type $G_2$, $F_4$, $E_6$ or
  $E_7$, then the centralizer of every semisimple element of $G$ cannot have factors of type $E_8$ (for dimensional
  reasons).

  Therefore, the previous proposition and the lemma in Section \ref{red-Z-section} together imply that
  $\widetilde{\tau}_G(h,\pi)\in\overline{\mathbb{Z}}[\frac{1}{M_G}]$ for all $h\in G^F$ and all
  $\pi \in \mathsf{K}^\circ_{G^\ast}$. We then deduce from (\ref{rest-tau-A}) that
 \begin{equation}\label{M-cor1} \mbox{$\tau_G(h\pi)=\widetilde{\tau}_G(h,\pi)\in\Lambda[\frac{1}{M_G}]$
      for all $h\in\Lambda\mathsf{E}_G$ and all $\pi \in \mathsf{K}^\circ_{G^\ast}$.}
  \end{equation}
  Now fit $G$ into the exact sequence (\ref{z-ext-D}). As $H$ therein has the same type of root datum as $G$, we have
  $M_H=M_G$, so (\ref{M-cor1}) applied to $H$ gives 
  $\tau_H(h\pi)\in\Lambda[\frac{1}{M_G}]$ for all
  $h\in\Lambda\mathsf{E}_H$ and all $\pi\in \mathsf{K}^\circ_{H^\ast} = \mathsf{K}_{H^\ast}$. 
  For our $G$, (\ref{pi-lambda-red}) then tells us that
  (\ref{key-lemma}) is true when $\Lambda$ therein is replaced by $\Lambda[\frac{1}{M_G}]$. Consequently, when all bad
  prime numbers for $G$ are invertible in $\Lambda$, we have $\Lambda[\frac{1}{M_G}]=\Lambda$ and
  $\Lambda\mathsf{E}_G=\Lambda\mathsf{K}_{G^\ast}$. \qed

\end{context}

{\small Tzu-Jan Li: Institute of Mathematics, Academia Sinica, 6F, Astronomy-Mathematics Building, No.\;1, Sec.\;4, Roosevelt Road, Taipei 106319, TAIWAN. \\ \indent {\it Email address}: {\tt tzujan@gate.sinica.edu.tw}

Jack Shotton: Department of Mathematical Sciences, Mathematical Sciences \& Computer Science Building, Durham University
Upper Mountjoy Campus, Stockton Road, Durham DH1 3LE UNITED KINGDOM\\ \indent {\it Email address}: {\tt
jack.g.shotton@durham.ac.uk}

}


\begin{thebibliography}{9}
\bibitem[BoKe]{Bonnafe-Kessar} C.\;Bonnaf{\'e} and R.\;Kessar, {\it On the endomorphism algebras of modular
Gelfand--Graev representations}, J.\;Alg.\;320 (2008)
\bibitem[Ca]{Carter} R.\;W.\;Carter, {\it On the Representation Theory of the Finite Groups of Lie Type over an
Algebraically Closed Field of Characteristic 0}, Encyclopaedia of Mathematical Sciences Vol.\;77, Algebra IX, Springer
(1996)
\bibitem[Cu]{Curtis} C.\;W.\;Curtis, {\it On the Gelfand--Graev Representations of a Reductive Group over a Finite
Field}, J.\;Alg.\;157 (1993)
\bibitem[De]{Deligne}
P.\;Deligne, {\it La conjecture de Weil : I},
Publ. math. I.H.{\'E}.S. 43 (1974)
\bibitem[DeLu]{Deligne--Lusztig} P.\;Deligne and G.\;Lusztig, {\it Representations of reductive groups over finite
fields}, Ann.\;Math.\;103-1 (1976)
\bibitem[DiMi]{Digne-Michel} 
F.\;Digne and J.\;Michel, 
Representations of finite groups of Lie Type, 2nd.\;ed., London
Math.\;Soc.\;Student Texts 95, Cambridge Univ.\;Press (2020)
\bibitem[DLM]{Digne-Lehrer-Michel} F.\;Digne, G.\;I.\;Lehrer and J.\;Michel, {\it The characters of the group of
rational points of a reductive group with non-connected centre}, J.\;reine angew.\;Math.\;425 (1992)
\bibitem[Hel]{Helm} D.\;Helm, {\it Curtis homomorphisms and the integral Bernstein center for $\mathrm{GL}_n$}, Alg.\;\&
Num.\;Th.\;14-10 (2020)
\bibitem[Her]{Herzig} F.\;Herzig, {\it The weight in a Serre-type conjecture for tame n-dimensional Galois
representations}, Duke Math.\;J.\;149-1 (2009)
\bibitem[Ja]{Jantzen} J.\;C.\;Jantzen, Representations of Algebraic Groups, Academic Press.\;Inc.\;(1987)

\bibitem[Li]{Thesis} T.-J.\;Li, {\it Sur l'alg{\`e}bre d'endomorphismes des repr{\'e}sentations de Gelfand--Graev
et le $\ell$-bloc unipotent de $\mathrm{GL}_2$ $p$-adique avec $\ell\neq p$}, Th{\`e}se de doctorat, Sorbonne Univ.\;(2022)
\bibitem[Li2]{Endomorphism} T.-J.\;Li, {\it On endomorphism algebras of Gelfand--Graev representations}, Represent.\;Theory 27 (2023)

\bibitem[Lu]{Lusztig} G.\;Lusztig, {\it Characters of {{Reductive Groups}} over a {{Finite Field}}.}, Ann. Math Studies
107 (1984)
\bibitem[HeMo]{Helm-Moss} D.\;Helm and G.\;Moss, {\it Converse theorems and the local Langlands correspondence in
families}, Inventiones Math. 214 (2018)
\bibitem[Se]{Serre} J.-P.\;Serre, {\it Linear Representations of Finite Groups}, Graduate Texts in Mathematics 42, Springer (1977)
\bibitem[Sh]{Shoji} T.\;Shoji, {\it Character Sheaves and Almost Characters of Reductive Groups, II}, Adv. Math. 111
(1995)
\bibitem[Sp]{Springer66} T.\;A.\;Springer, {\it Some Arithmetical Results on Semi-Simple Lie Algebras},
Publ. math. I.H.\'E.S. 30 (1966)
\bibitem[Sp2]{Springer} T.\;A.\;Springer, {\it A Construction of Representations of Weyl Groups}, Inventiones Math.\;44
(1978)
\bibitem[St]{Steinberg} R.\;Steinberg, {\it Representations of algebraic groups}, J.\;Nagoya Math.\;22 (1963)


\end{thebibliography}
\end{document}